UDC 519.652 + 517.518.85

# Trigonometric splines and some of their properties

Denysiuk V.P. Dr of phys-math. sciences, Professor, Kiev, Ukraine

National Aviation University

kvomden@nau.edu.ua

**Annotation**

A class of trigonometric interpolation splines depending on parameter vectors, selected convergence factors and interpolation factors is considered. The concept of crosslink grids and interpolation grids is introduced; these grids can match or y not match. Parseval equality is given and some properties of splines are discussed. The material is illustrated by a numerical example.

**Keywords**:

Generalized trigonometric functions, interpolation, polynomial and trigonometric splines.

**Introduction**

Approximation, respectively, the representation of a known or unknown function through a set of some special functions can be considered as a central topic of analysis; such special functions are well defined, easy to calculate, and have certain analytical properties [ 1]. Algebraic and trigonometric polynomials, exponential functions [2], polynomial [3] and trigonometric [4] splines, etc. often act as special functions.

In applications, time often acts as an argument for the functions under study $t$, $t \in [0,T]$; according to it, we will denote the functions in the future $f(t)$.

In many cases, the function under study $f(t)$, which, as a rule, is at least continuous, has a rather complex analytical representation or such representation is unknown at all. In such cases, the function is discretized, namely replaced by a finite sequence of its instantaneous values at the nodes of some, as a rule, equidistant grid specified on $[0,T]$. In the future there is a problem of approximation of this sequence in one sense or another by a linear combination of special functions of a certain class. In this approach, it is important to choose a class of special functions, which minimizes the deviation of the function $f(t)$ from an approximating function built on discrete data.

It is clear that without imposing restrictions on the class of functions $f(t)$, it is difficult to make such a choice. We will consider the case when functions $f(t) \in W_v^r$, where $W_v^r$ - a class of periodic functions having an absolutely continuous derivative of order $r-1$ ($r = 1, 2, ...$), and a derivative of the order $r$ is a function of limited variation.

It is known [5,6] that the best apparatus for approximating the functions of a class $W_v^r$ there are simple polynomial splines. The theory of such splines is well developed (see, for example, [3], [7], [8], [9]), etc. The main disadvantage of polynomial splines, in our opinion, is that they have a lump structure; this leads to the fact that in practice mainly third-order splines are used, which are sewn together from pieces of algebraic polynomials of the third degree. However, this structure of splines significantly limits their application in many problems of computational mathematics.

By Schoenberg [10] were proposed trigonometric splines, which are crosslinked from trigonometric polynomials of a certain degree, and have the same disadvantage as polynomial splines - a lump structure.

In [4], [11], [12] another principle of construction of trigonometric splines was proposed, their representation by uniformly convergent trigonometric series (Fourier series), the coefficients of which have a certain order of decrease. The undeniable advantage of such splines is that they are given as a single expression over the entire interval of the function. One of the properties of the trigonometric splines introduced in this way is the certain values of the parameters they coincide with simple polynomial splines [11]; therefore, in this case, all the results of approximation estimates obtained for polynomial splines can be transferred to trigonometric splines. Moreover, the well-developed theory of Fourier trigonometric series gives reason to expect the appearance of new results in this area.

## The purpose of the work

Construction of a class of interpolation trigonometric splines, which depend on the vectors of parameters $\Gamma$ and $H$, and study of some of their properties.

## The main part

In [4] interpolation trigonometric splines were considered

$$St^{(I_1,I_2)}(\nu,r,N,t) = \frac{a_0^{(I_2)}}{2} + \sum_{k=1}^{\frac{N-1}{2}} \frac{1}{h^{(I_1,I_2)}(r,k)} \left[ a_\kappa^{(I_2)} C_k^{(I_1)}(\nu,r,N,t) + b_\kappa^{(I_2)} S_k^{(I_1)}(\nu,r,N,t) \right]$$

where

$$C_k^{(I_1)}(\nu,r,N,t) = \nu_k(r)\cos kt + \sum_{m=1}^{\infty}(-1)^{mI_1}\left[\nu_{mN+k}(r)\cos(mN+k)t + \nu_{mN-k}(r)\cos(mN-k)t\right];$$

$$S_k^{(I_1)}(\nu,r,N,t) = \nu_k(r)\sin kt + \sum_{m=1}^{\infty}(-1)^{mI_1}\left[\nu_{mN+k}(r)\sin(mN+k)t - \nu_{mN-k}(r)\sin(mN-k)t\right], \quad (1)$$

convergence factors $\nu_k(r)$, having a descending order $O(k^{-(1+r)})$, and an interpolation factor

$$h^{(I_1,I_2)}(r,k) = \nu_k(r) + \sum_{m=1}^{\infty}(-1)^{m(I_1-2I_1I_2+I_2)}\left[\nu_{mN+k}(r) + \nu_{mN-k}(r)\right], \quad (2)$$

which provide interpolation of setpoints in the grid nodes $\Delta_N^{(I)} = \{t_i^{(I)}\}_{i=1}^N$, ($I=0,1$), where $t_j^{(0)} = \frac{2\pi}{N}(j-1)$, but $t_j^{(1)} = \frac{\pi}{N}(2j-1)$, and $N=2n+1$, ($n=1,2,...$), and for which their polynomial analogues exist (or can be built).

However, it is easy to construct other types of trigonometric interpolation splines, which may be of interest in certain interpolation problems.

We consider two vectors $\Gamma = \{\gamma_1,\gamma_2,\gamma_3\}$ and $H = \{\eta_1,\eta_2,\eta_3\}$ and consider the spline

$$St^{(I_1,I_2)}(\Gamma,H,\nu,r,N,t) = \frac{a_0^{(I_2)}}{2} + \sum_{k=1}^{\frac{N-1}{2}} \left[ a_\kappa^{(I_2)} \frac{C_k^{(I_1)}(\Gamma,\nu,r,N,t)}{hc^{(I_1,I_2)}(\Gamma,r,k)} + b_\kappa^{(I_2)} \frac{S_k^{(I_1)}(H,\nu,r,N,t)}{hs^{(I_1,I_2)}(H,r,k)} \right]$$

functions

$$C_k^{(I_1)}(\Gamma,\nu,r,N,t) =$$
$$= \gamma_1\nu_k(r)\cos kt + \sum_{m=1}^{\infty}(-1)^{mI_1}\left[\gamma_2\nu_{mN+k}(r)\cos(mN+k)t + \gamma_3\nu_{mN-k}(r)\cos(mN-k)t\right];$$

$$S_k^{(I_1)}(H,\nu,r,N,t) = \quad (3)$$
$$= \eta_1\nu_k(r)\sin kt + \sum_{m=1}^{\infty}(-1)^{mI_1}\left[\eta_2\nu_{mN+k}(r)\sin(mN+k)t - \eta_3\nu_{mN-k}(r)\sin(mN-k)t\right],$$

with convergence factors $\nu_k(r)$, having a descending order $O(k^{-(1+r)})$, and interpolation factors

$$hc^{(I_1,I_2)}(\Gamma,\nu,r,N,k) = \gamma_1\nu_k(r) + \sum_{m=1}^{\infty}(-1)^{m(I_1-2I_1I_2+I_2)}\left[\gamma_2\nu_{mN+k}(r) + \gamma_3\nu_{mN-k}(r)\right],$$

$$hs^{(I_1,I_2)}(H,\nu,r,N,k) = \eta_1\nu_k(r) + \sum_{m=1}^{\infty}(-1)^{m(I_1-2I_1I_2+I_2)}\left[\eta_2\nu_{mN+k}(r) + \eta_3\nu_{mN-k}(r)\right]. \quad (4)$$

where parameters $\gamma_k$ i $\eta_k$, ($k=1,2,3$) take arbitrary real values, and these parameters don't turn simultaneously to 0. Note that in this paper we do not consider the case when the parameters $\eta_k$ i $\theta_k$, ($k=1,2,3$) take complex values, leaving the consideration of this case for later.

It is clear that functions (3) and interpolation factors (4) differ from (1) and (2) by the presence of parameters $\gamma_k$ i $\eta_k$, ($k=1,2,3$).

If $\gamma_k = 1$ and $\eta_k = 1$, ($k = 1, 2, 3$), then functions (3) will be called simple basic functions, and trigonometric splines constructed using simple basic functions with corresponding interpolation factors (4) will be called simple trigonometric splines. Simultaneously, time functions $C_k^{(I)}(\Gamma, \nu, r, N, t)$, and $S_k^{(I)}(H, \nu, r, N, t)$ and their corresponding multipliers $hc^{(I)}(\Gamma, r, k)$ and $hs^{(I)}(H, r, k)$, we will denote accordingly $C_k^{(I)}(\nu, r, N, t)$, $S_k^{(I)}(\nu, r, N, t)$ and $hc^{(I)}(r, k)$ i $hs^{(I)}(r, k)$. It is clear that in this case functions (3) and interpolation factors (4) coincide with functions (1) and factors (2).

Note that in the case where $\gamma_1 = \eta_1 = 1$, a , we have a common interpolation trigonometric polynomial.

As follows from (3), (4), in the general case, trigonometric splines depend from the parameters $\gamma_k$ and $\eta_k$, ($k = 1, 2, 3$), type of selected convergence factors $\nu_k(r)$, crosslinking grid and interpolation grid. In the future, the index of the crosslink grid will be denoted as $I_1$, and the index of interpolation grid as $I_2$, ($I_1, I_2 = 0, 1$). Note that the term crosslink grid refers to the grid on which polynomial analogues (if they exist) of trigonometric splines are crosslinked. It is clear that the crosslink grid is determined by functions $C_k^{(I_1)}(\Gamma, \nu, r, N, t)$ and $S_k^{(I_1)}(H, \nu, r, N, t)$; the interpolation grid is determined by multipliers $hc^{(I_1, I_2)}(\Gamma, r, k)$ and $hs^{(I_1, I_2)}(H, r, k)$.

Note that the crosslink grid and the interpolation grid may match or not match. Thus, the construction of trigonometric splines consists of the following stages.

1. Set the system of parameters $\gamma_k$ and $\eta_k$, ($k = 1, 2, 3$).
2. Choose the type of convergence factors $\nu_k(r, N)$; for the role of such factors can be selected, for example, type factors $\psi(r, k) / k^{1+r}$, where $\psi(r, k)$ - some limited function of its arguments, as well as Fourier coefficients of some functions of classes $W_\nu^r$ [12], which have a descending order $O(k^{-(1+r)})$.
3. Choose the crosslink grid, which mean - set the parameter $I_1$;
4. Select the interpolation grid, which mean set the parameter $I_2$;
5. Determine the number of continuous derivatives of the trigonometric spline, is mean set the parameter $r$ ($r = 0, 1, 2, \ldots$); note that in the case $r = 0$ the convergence of trigonometric series given by trigonometric splines is not even; however, it is often advisable to consider this case as well.

Consider an example that illustrates the approach of construction of trigonometric splines in more detail.

It is clear that to build trigonometric interpolation splines, you must first specify the number of nodes $N$, ($N = 2n + 1$), ($n = 1, 2, \ldots$), on the segment $[0, 2\pi)$ and thereby determine the crosslink grids and interpolation grids. You must also specify the value of the interpolated function on these grids. So, let it be in $[0, 2\pi)$ given grids: $\Delta_N^{(I)} = \{t_j^{(I)}\}$, ($t_j^{(0)} = \frac{2\pi}{N}(j-1)$, $t_j^{(1)} = \frac{\pi}{N}(2j-1)$). Suppose that in this interval some periodic continuous function is also given $f(t)$ and let the values be known $\{f(t_j^{(I)})\}_{j=1}^N = \{f_j^{(I)}\}_{j=1}^N$ of this function in grid nodes $\Delta_N^{(I)}$. As you know, the coefficients $a_0^{(I)}, a_k^{(I)}, b_k^{(I)}$, ($k = 1, 2, \ldots, n$) trigonometric polynomial

$$T_n^{(I)}(t) = \frac{a_0^{(I)}}{2} + \sum_{k=1}^{\frac{N-1}{2}} \left(a_k^{(I)} \cos kt + b_k^{(I)} \sin kt\right)$$

interpolating function $f(t)$ in grid nodes $\Delta_N^{(I)}$, calculated by formulas

$$a_0^{(I)} = \frac{1}{N} \sum_{j=1}^N f_j^{(I)};$$

$$a_k^{(I)} = \frac{2}{N} \sum_{j=1}^N f_j^{(I)} \cos kt_j^{(I)}; \quad b_k^{(I)} = \frac{2}{N} \sum_{j=1}^N f_j^{(I)} \sin kt_j^{(I)} \quad (6)$$

For certainty we will put $N = 9$, and $\{f_j^{(I)}\}_{j=1}^{9} = \{3,1,3,2,4,1,3,1,2\}$.

We now turn to the definition of the type and parameters of the trigonometric spline.
1. Put for certainty $\gamma_k = 1$ and $\eta_k = 1$, ($k = 1,2,3$).
2. Choose the convergence factor

$$v_k(r,N) = \left[\frac{\sin\left(\frac{\alpha}{2}k\right)}{\frac{\alpha}{2}k}\right]^{1+r}, \qquad (7)$$

where $\alpha = 2\pi/N$; it is clear that the order of decline of this multiplier $O(k^{-(1+r)})$. Note that the role of the parameter $\alpha$ requires a separate study and is not considered here.

3. Select the crosslink grid - set the parameter $I_1$, ($I_1 = 0,1$); later we will consider both values of this parameter.
4. Select the interpolation grid - set the parameter $I_2$, ($I_2 = 0,1$); later we will consider both values of this parameter.
5. Set the parameter $r$ ($r = 0,1,2,...$); later we will consider some values of this parameter.
6. According to the proposed classification, we consider simple trigonometric splines, which are denoted as follows

$$St^{(I_1,I_2)}(v,r,N,t),$$

Where index $I_1$ determines the crosslink grid, and the index $I_2$ determines the interpolation grid $(I_1, I_2 = 0,1)$, parameter $v$ determines the type of selected convergence factors, parameter $r$ ($r = 1,2,...$) determines the smoothness of the spline, the parameter $N$ ($N = 2n+1$, $n = 1,2,...$) determines the number of nodes of the uniform interpolation grid specified on the interval $[0, 2\pi)$, and $t$ - spline argument.

Here are the graphs of these splines for some parameter values and different values of the grid indexes $I_1$, $I_2$, taking into account the following remarks.

1. Graphs 1-4 show the grid nodes with vertical lines $\Delta_N^{(0)}$.
2. Grid nodes $\Delta_N^{(1)}$ located between the nodes of the grid $\Delta_N^{(0)}$ and graphs 1-4 are not displayed..
3. Graphs 1-4 instead of complex spline symbols $St^{(I_1,I_2)}(v,r,N,t)$ simplified type designations are used $StI_1I_2(r,t)$, corresponding to different combinations of index values $I_1, I_2$.

Let's move on to consider the trigonometric spline $St^{(I_1,I_2)}(v,r,N,t)$.

It is clear that when $I_1 = 0$, $I_2 = 0$ we get;

$$C_k^{(0)}(v,r,N,t) = v_k(r)\cos kt + \sum_{m=1}^{\infty}\left[v_{mN+k}(r)\cos(mN+k)t + v_{mN-k}(r)\cos(mN-k)t\right];$$

$$S_k^{(0)}(v,r,N,t) = v_k(r)\sin kt + \sum_{m=1}^{\infty}\left[v_{mN+k}(r)\sin(mN+k)t - v_{mN-k}(r)\sin(mN-k)t\right],$$

$$hc^{(0,0)}(r,k) = v_k(r) + \sum_{m=1}^{\infty}\left[v_{mN+k}(r) + v_{mN-k}(r)\right]$$

$$hs^{(0,0)}(r,k) = v_k(r) + \sum_{m=1}^{\infty}\left[v_{mN+k}(r) + v_{mN-k}(r)\right]$$

and a trigonometric spline $St^{(0,0)}(v,r,N,t)$

$$St^{(0,0)}(v,r,N,t) = \frac{a_0^{(0)}}{2} + \sum_{k=1}^{\frac{N-1}{2}}\left[a_k^{(0)}\frac{C_k^{(0)}(v,r,N,t)}{hc^{(0,0)}(r,k)} + b_k^{(0)}\frac{S_k^{(0)}(v,r,N,t)}{hs^{(0,0)}(r,k)}\right].$$

Graphs of this spline for some parameter $r$ are shown in Fig.1

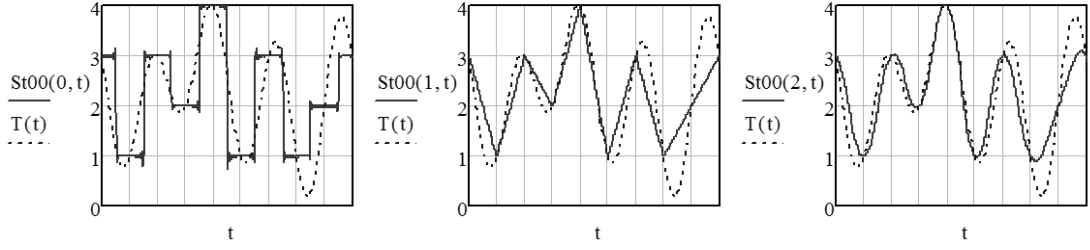

Fig.1 Graphs of the spline $St^{(0,0)}(v,r,N,t)$ for $r=0,1,3$. Crosslinking and interpolation carried out in the nodes of the grid $\Delta_N^{(0)}$.

For the spline $St^{(0,1)}(v,r,N,t)$ while $I_1=0$, $I_2=1$, we get;

$$C_k^{(0)}(v,r,N,t) = v_k(r)\cos kt + \sum_{m=1}^{\infty}\left[v_{mN+k}(r)\cos(mN+k)t + v_{mN-k}(r)\cos(mN-k)t\right];$$

$$S_k^{(0)}(v,r,N,t) = v_k(r)\sin kt + \sum_{m=1}^{\infty}\left[v_{mN+k}(r)\sin(mN+k)t - v_{mN-k}(r)\sin(mN-k)t\right],$$

$$hc^{(0,1)}(r,k) = v_k(r) + \sum_{m=1}^{\infty}(-1)^m\left[v_{mN+k}(r) + v_{mN-k}(r)\right],$$

$$hs^{(0,1)}(r,k) = v_k(r) + \sum_{m=1}^{\infty}(-1)^m\left[v_{mN+k}(r) + v_{mN-k}(r)\right],$$

and a trigonometric spline $St^{(0,1)}(v,r,N,t)$

$$St^{(0,1)}(v,r,N,t) = \frac{a_0^{(1)}}{2} + \sum_{k=1}^{\frac{N-1}{2}}\left[a_k^{(1)}\frac{C_k^{(0)}(v,r,N,t)}{hc^{(0,1)}(r,k)} + b_k^{(1)}\frac{S_k^{(0)}(v,r,N,t)}{hs^{(0,1)}(r,k)}\right].$$

Graphs of this trigonometric spline are shown in Fig. 2.

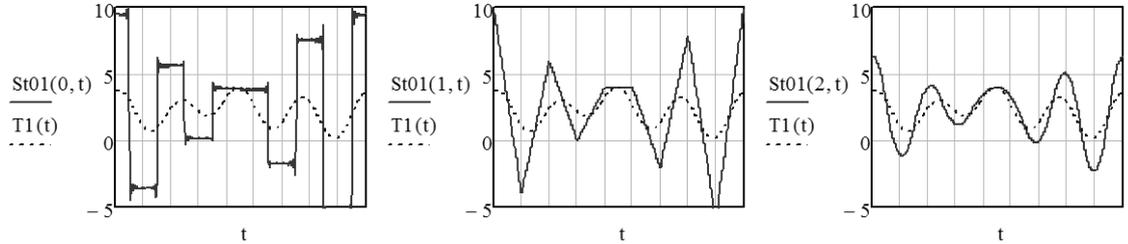

Fig.2 Graphs of the spline $St^{(0,1)}(v,r,N,t)$ for $r=0,1,3$. Crosslinking carried out in the nodes of the grid $\Delta_N^{(0)}$, and interpolation at grid nodes $\Delta_N^{(1)}$

For the spline $St^{(1,0)}(v,r,N,t)$ we get;

$$C_k^{(1)}(v,r,N,t) = v_k(r)\cos kt + \sum_{m=1}^{\infty}(-1)^m\left[v_{mN+k}(r)\cos(mN+k)t + v_{mN-k}(r)\cos(mN-k)t\right];$$

$$S_k^{(1)}(v,r,N,t) = v_k(r)\sin kt + \sum_{m=1}^{\infty}(-1)^m\left[v_{mN+k}(r)\sin(mN+k)t - v_{mN-k}(r)\sin(mN-k)t\right],$$

$$hc^{(1,0)}(v,r,N,k) = v_k(r) + \sum_{m=1}^{\infty}(-1)^m\left[v_{mN+k}(r) + v_{mN-k}(r)\right],$$

$$hs^{(1,0)}(v,r,N,k) = v_k(r) + \sum_{m=1}^{\infty}(-1)^m\left[v_{mN+k}(r) + v_{mN-k}(r)\right],$$

and a trigonometric spline $St^{(1,0)}(v,r,N,t)$

$$St^{(1,0)}(\nu,r,N,t) = \frac{a_0^{(0)}}{2} + \sum_{k=1}^{\frac{N-1}{2}} \left[ a_k^{(0)} \frac{C_k^{(1)}(\nu,r,N,t)}{hc^{(1,0)}(\nu,r,N,k)} + b_k^{(0)} \frac{S_k^{(1)}(\nu,r,N,t)}{hs^{(1,0)}(\nu,r,N,k)} \right].$$

Graphs of this trigonometric spline are shown in Fig. 3.

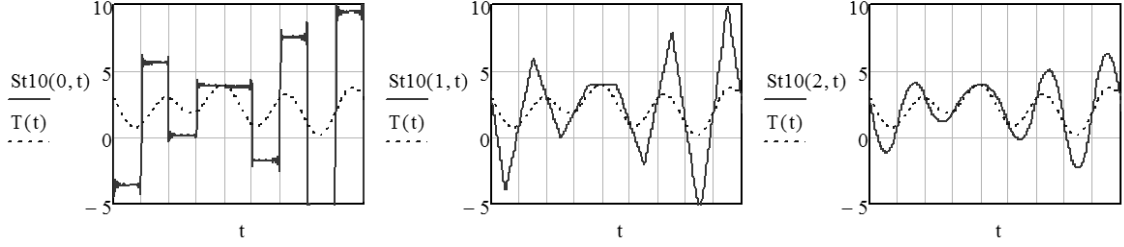

Fig.3 Graphs of the spline $St^{(1,0)}(\nu,r,N,t)$ for $r=0,1,3$. Crosslinking carried out in the nodes of the grid $\Delta_N^{(1)}$, and interpolation at grid nodes $\Delta_N^{(0)}$ /

For the spline $St^{(1,1)}(\nu,r,N,t)$ we get;

$$C_k^{(1)}(\nu,r,N,t) = \nu_k(r)\cos kt + \sum_{m=1}^{\infty} (-1)^m \left[ \nu_{mN+k}(r)\cos(mN+k)t + \nu_{mN-k}(r)\cos(mN-k)t \right];$$

$$S_k^{(1)}(\nu,r,N,t) = \nu_k(r)\sin kt + \sum_{m=1}^{\infty} (-1)^m \left[ \nu_{mN+k}(r)\sin(mN+k)t - \nu_{mN-k}(r)\sin(mN-k)t \right],$$

$$hc^{(1,1)}(r,k) = \nu_k(r) + \sum_{m=1}^{\infty} \left[ \nu_{mN+k}(r) + \nu_{mN-k}(r) \right],$$

$$hs^{(1,1)}(r,k) = \nu_k(r) + \sum_{m=1}^{\infty} \left[ \nu_{mN+k}(r) + \nu_{mN-k}(r) \right],$$

and a trigonometric spline $St^{(1,1)}(\nu,r,N,t)$

$$St^{(1,1)}(\nu,r,N,t) = \frac{a_0^{(1)}}{2} + \sum_{k=1}^{\frac{N-1}{2}} \left[ a_k^{(1)} \frac{C_k^{(1)}(\nu,r,N,t)}{hc^{(1,1)}(r,k)} + b_k^{(1)} \frac{S_k^{(1)}(\nu,r,N,t)}{hs^{(1,1)}(r,k)} \right].$$

Graphs of this trigonometric spline are shown in Fig.4.

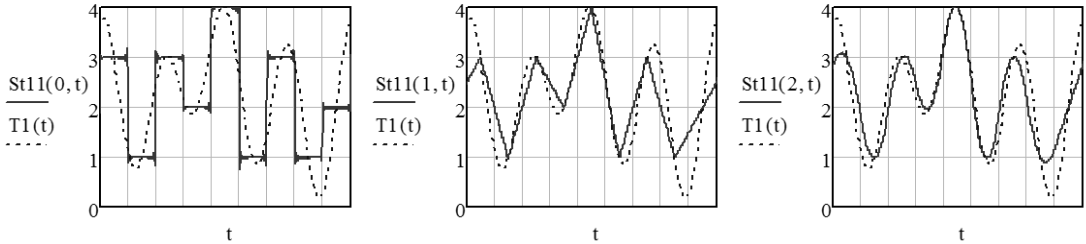

Fig.4 Spline graphs $St^{(1,1)}(\nu,r,N,t)$ for $r=0,1,3$. Crosslinking and interpolation carried out in the nodes of the grid $\Delta_N^{(1)}$.

Thus, giving different values to the indexes $I_1$ and $I_2$, you can choose crosslinking and interpolation grids based on the needs of the problem.

Consider now the case where the parameters $\gamma_k$ i $\eta_k$, ($k=1,2,3$) not equal to 1. Let's consider, for example, $\gamma_1 = -.5$, $\gamma_2 = 1.5$, $\gamma_3 = -.7$, $\eta_1 = .3$, $\eta_2 = -.7$, $\eta_3 = -1.5$ and give graphs of trigonometric splines $St^{(I,I)}(\Gamma,H,\nu,r,N,t)$; we will use abbreviated symbols on the graphs, detailing them on the signatures to the graphs.

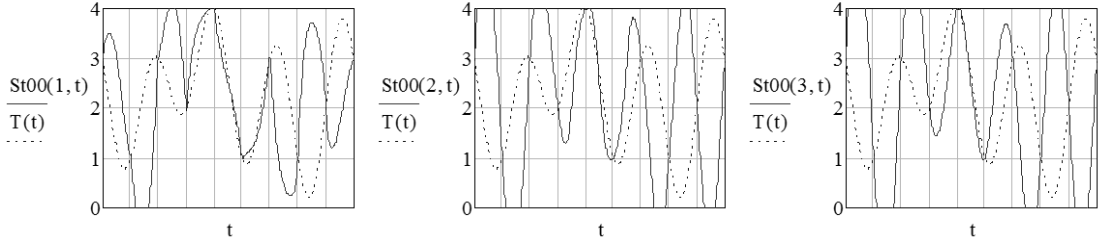

Fig.5. Trigonometric spline graphs $St^{(0,0)}(\Gamma, H, \nu, r, 9, t)$ for $r = 1, 2, 3$; crosslinking and interpolation in grid nodes $\Delta_N^{(0)}$.

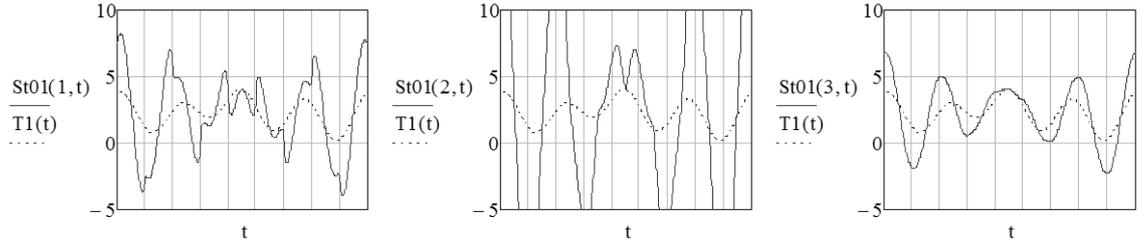

Fig.6. Trigonometric spline graphs $St^{(0,1)}(\Gamma, H, \nu, r, 9, t)$ for $r = 1, 2, 3$; crosslinking in grid nodes $\Delta_N^{(0)}$, and interpolation at grid nodes $\Delta_N^{(1)}$.

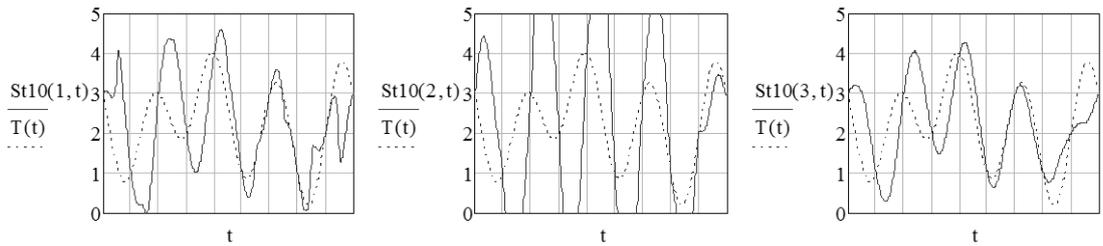

Fig.7. Trigonometric spline graphs $St^{(1,0)}(\Gamma, H, \nu, r, 9, t)$ for $r = 1, 2, 3$; crosslinking in grid nodes $\Delta_N^{(1)}$, and interpolation at grid nodes $\Delta_N^{(0)}$.

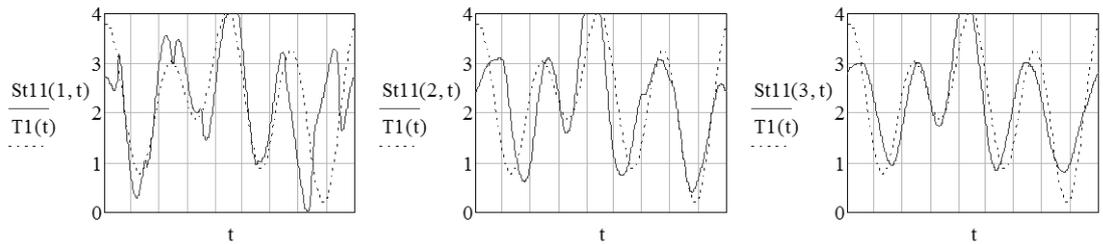

Fig.8. Trigonometric spline graphs $St^{(1,0)}(\Gamma, H, \nu, r, 9, t)$ for $r = 1, 2, 3$; Crosslinking and interpolation carried out in the nodes of the grid $\Delta_N^{(1)}$.

Analysing the graphs in Fig.5-8, it is advisable to pay attention to the rather irregular behavior of trigonometric splines with such parameter values $\Gamma$ and $H$; however, our goal was only to illustrate the construction of trigonometric splines at selected values of these parameters.

Let's consider some properties of trigonometric interpolation splines.

In the theory of polynomial splines, half-norms of splines and their derivatives in space $L_{[a,b]}^p$ are often considered, which is defined as follows

$$\|f\|_p = \left\{\int_a^b |f(x)|^p \right\}^{\frac{1}{p}} dx \ , \ (p = 1, 2, \ldots). \tag{8}$$

The most popular case is $p = 2$, because, firstly, in this case, the half-norm has a simple interpretation, and secondly, such a half-norm is generated by a scalar product. Thus, in particular, in [7] the extreme properties of polynomial simple interpolation splines of odd order are considered when $p = 2$

Consider some properties of trigonometric interpolation splines depending on their degree $r$ ($r = 1, 2, \ldots$), interpolation grids and crosslink grids and some vectors $\Gamma$ and $H$. ;

First of all, we note that in the case of trigonometric interpolation splines it is expedient to consider the norm

$$\|f\| = \frac{1}{\pi} \int_0^{2\pi} |f(x)|^2 dx \tag{9}$$

which differs from the half-norm (8) by the presence of a constant factor and the absence of an exponent .5. The expediency of using such a norm is explained by the fact that many results in the theory of trigonometric series are formulated in terms of such a norm; examples include Bessel's inequality, Parseval's equality, etc. Note that (9) is often interpreted as the average power of the function $f(x)$ on the period $[0, 2\pi)$.

It is clear that for each trigonometric spline it is easy to calculate the value

$$PSt^{(I_1,I_2)}(\Gamma,H,\nu,r,N) = \frac{1}{\pi} \int_0^{2\pi} \left|St^{(I_1,I_2)}(\Gamma,H,\nu,r,N,t)\right|^2 dt,$$

As trigonometric splines at $r \geq 1$ are represented by uniformly convergent trigonometric series, then these series must be Fourier series of these splines; so, there is a Parseval equality that can be represented as

$$\frac{1}{\pi} \int_0^{2\pi} \left|St^{(I_1,I_2)}(\Gamma,H,\nu,r,N,t)\right|^2 dt =$$

$$= \frac{\left[a_0^{(I_2)}\right]^2}{2} + \sum_{k=1}^{\frac{N-1}{2}} \left[\left[a_k^{(I_2)}\right]^2 PC_k^{(I_1,I_2)}(\Gamma,\nu,r,N,t) + \left[b_k^{(I_2)}\right]^2 PS_k^{(I_1,I_2)}(H,\nu,r,N,t)\right],$$

where

$$PC_k^{(I_1,I_2)}(\Gamma,\nu,r,N,t) =$$
$$= \frac{1}{\left[HC^{(I_1,I_2)}(\Gamma,\nu,r,N,k)\right]^2} \left\{\left[\gamma_1 \nu_k(r)\right]^2 + \sum_{m=1}^{\infty}(-1)^{mI}\left[\left[\gamma_2 \nu_{mN+k}(r)\right]^2 + \left[\gamma_3 \nu_{mN-k}(r)\right]^2\right]\right\};$$

$$PS_k^{(I_1,I_2)}(H,\nu,r,N,t) =$$
$$= \frac{1}{\left[HS^{(I_1,I_2)}(H,\nu,r,N,k)\right]^2} \left\{\left[\eta_1 \nu_k(r)\right]^2 + \sum_{m=1}^{\infty}(-1)^{mI}\left[\left[\eta_2 \nu_{mN+k}(r)\right]^2 + \left[\eta_3 \nu_{mN-k}(r)\right]^2\right]\right\}.$$

First of all, note that when vectors $\Gamma$ and $H$ are equal, then there are such equalities
$$PSt^{(0,0)}(\Gamma,H,\nu,r,N) = PSt^{(1,1)}(\Gamma,H,\nu,r,N);$$
$$PSt^{(0,1)}(\Gamma,H,\nu,r,N) = PSt^{(1,0)}(\Gamma,H,\nu,r,N).$$

Consider the magnitude $PSt^{(I_1,I_2)}(\Gamma,H,\nu,r,N)$ taking to consideration $\nu_k(r,N)$ and unit vectors $\Gamma$ and $H$; as we have already seen, in this case at odd values of the parameter $r$ and $I_2 = 0$ we have simple polynomial splines of odd degree, and for even values of this parameter and $I_2 = 1$ - polynomial splines of even degree.

Due to the fact that [13]

$$\lim_{r \to \infty} St^{(0,I_2)}(\Gamma, H, \nu, r, N, t) = \lim_{r \to \infty} St^{(1,I_2)}(\Gamma, H, \nu, r, N, t) = T_N^{I_2}(t)$$

it follows that

$$\lim_{r \to \infty} PSt^{(0,I_2)}(\Gamma, H, \nu, r, N, t) = \lim_{r \to \infty} PSt^{(1,I_2)}(\Gamma, H, \nu, r, N, t) = PT_N^{I_2}(t),$$

and values $PSt^{(0,0)}(\Gamma, H, \nu, r, N)$ and $PSt^{(1,1)}(\Gamma, H, \nu, r, N)$ go to this boundary from below, and the magnitude $PSt^{(0,1)}(\Gamma, H, \nu, r, N)$ and $PSt^{(1,0)}(\Gamma, H, \nu, r, N)$ - from above.

Changing the components of vectors $\Gamma$ and $H$ you can get other results. Thus, in particular, when $\Gamma = \{.5, 10, 1\}$ and $\Gamma = H$, splines were obtained whose average power was less than the average power of simple polynomial splines of odd order.

It is clear that the presence of parameter vectors $\Gamma$ and $H$ provides opportunities for setting optimization problems for these and other functionalities from trigonometric splines and their derivatives.

Finally, we note that in many cases trigonometric splines, which are periodic functions, have to approximate non-periodic functions. In this case, it is advisable to apply special methods of periodic continuation of non-periodic functions [14].

## Conclusions

1. The class of interpolation trigonometric splines depending on the parameter vectors $\Gamma$ and $H$ is considered.
2. A set of simple trigonometric splines obtained for single vectors $\Gamma$ and $H$ is introduced.
3. Simple trigonometric splines $St^{(0,0)}(\nu, r, N, t)$ at the convergence factor (7) and at odd values of the parameter $r$ have polynomial analogues - simple polynomial splines of odd degree.
3. Simple trigonometric splines $St^{(0,1)}(\nu, r, N, t)$ for the convergence factor (7) and for even values of the parameter, they have polynomial analogues — simple polynomial splines of even degree.
4. Other considered trigonometric splines have no polynomial analogues.
5. There are other types of convergence factors exists, the introduction of which into simple trigonometric splines leads to simple polynomial splines..
6. Since interpolation trigonometric splines are determined by the coefficients of trigonometric interpolation polynomials, known FFT (fast Fourier transform) algorithms can be widely used in their construction.
6. The presence of parameter vectors $\Gamma$ and $H$ makes it possible to set optimization problems for certain functionalities given in the class of trigonometric splines.
7. Of course, the constructed class of trigonometric splines requires further research.

## List of references